\documentclass{article}
\usepackage{latexsym}
\usepackage{amssymb}
\usepackage{amsfonts}                     
\usepackage{amsmath}
\usepackage[all]{xy}
\usepackage{color}
\usepackage{mathrsfs}
\usepackage[utf8]{inputenc}
\usepackage[T1]{fontenc}
\usepackage{amscd}
\usepackage{amsxtra} 
\usepackage{url}
\newtheorem{thm}{\bf Theorem}[section]
\newtheorem{cor}[thm]{\bf Corollary}
\newtheorem{lem}[thm]{\bf Lemma}
\newtheorem{prop}[thm]{\bf Proposition}
\newtheorem{defn}[thm]{\bf Definition}

\newtheorem{exmp}[thm]{\bf Example}

\newcommand{\field}[1]{\mathbb{#1}}

\newcommand{\Q }{\field{Q}}
\newcommand{\Z }{\field{Z}}

\def\E{{\cal E}}

\def\X{{\cal X}}

\def\Ext{{\rm Ext}}
\def\Hom{{\rm Hom}}
\def\Im{{\rm Im}}

\def\proof{{\parindent0pt {\bf Proof.\ }}}
\def\Tor{{\rm Tor}}

\def\Ker{{\rm ker}}
\def\SFP{{\rm $S$-FP}}

\newcommand{\cqfd}
{\hspace{1cm}
\rule{2mm}{2mm}%
\medbreak%
\par%
}

\begin{document}

\title{On $S$-injective modules}

\author{Driss Bennis and Ayoub Bouziri}

\date{}
 
\maketitle
\begin{abstract} 
 Let $R$ be a commutative ring with identity, and let $S$ be a multiplicative subset of $R$. In this paper, we introduce the notion of $S$-injective modules as a weak version of injective modules. Among other results, we provide an $S$-version of Baer's characterization of injective modules. We also present an $S$-version of Lambek's characterization of flat modules: an $R$-module $M$ is $S$-flat if and only if its character, $\Hom_{\mathbb{Z}}(M,\mathbb{Q}/\mathbb{Z})$, is an $S$-injective $R$-module. As applications, we establish, under certain conditions, $S$-counterparts of the Cartan–Eilenberg–Bass and Cheatham-Stone characterizations of Noetherian rings.

  \end{abstract}

\medskip 
 
{\scriptsize \textbf{Mathematics Subject Classification (2020)}: 13C11, 13C13, 13E99. }

 {\scriptsize \textbf{Key Words}: $S$-Noetherian rings, $S$-flat modules, $S$-injective modules.}

\section{Introduction}
\hskip .5cm  Throughout this paper, $R$ is a commutative ring with identity, all modules are unitary and $S$ is a multiplicative subset of $R$; that is, $1 \in S$ and $s_1s_2 \in S$ for any $s_1 \in S, s_2 \in S$. Unless explicitly stated otherwise, when we refer to a multiplicative subset $S$ of $R$, we implicitly assume that $0 \notin S$. This assumption will be used in the sequel without explicit mention. Let $M$ be an $R$-module. As usual, we use $M^+$ and $M_S$ to denote, respectively, the character module $\Hom_{\Z}(M, \Q/\Z)$ and the localization of $M$ at $S$. Recall that $M_S \cong M \otimes_R R_S$.

\medskip 
In the last years, the notion of $S$-property draw attention of several authors. This  notion  was introduced in 2002 by D. D. Anderson and Dumitrescu where they defined the notions of $S$-finite modules and $S$-Noetherian rings. Namely, an $R$-module $M$ is said to be $S$-finite if there exist a finitely generated submodule $N$ of $M$ and $s \in S$ such that $sM \subseteq N$. A commutative ring $R$ is said to be $S$-Noetherian if every ideal of $R$ is $S$-finite \cite[Definition 1]{And2}.

 In \cite{Ben1}, Bennis and El Hajoui investigated an $S$-version of finitely presented modules and coherent rings which are called, respectively, $S$-finitely presented modules and $S$-coherent rings. An $R$-module $M$ is said to be $S$-finitely presented if there exists an exact sequence of $R$-modules $0 \to K\to L\to  M \to 0$, where $L$ is a  finitely generated free $R$-module and $K$ is an $S$-finite R-module. A  commutative ring $R$ is called $S$-coherent, if every finitely generated ideal of $R$ is $S$-finitely presented. They showed that the $S$-coherent rings have a similar characterization to the classical one given by Chase for coherent rings  \cite[Theorem 3.8]{Cha1}. Subsequently, they asked whether there exists an $S$-version of Chase's theorem \cite[Theorem 2.1]{Cha1}. In other words, how to define an $S$-version of flatness that characterizes $S$-coherent rings similarly to the classical case?  This problem was solved by the  notion    of $S$-flat module in \cite{Qi1}. Recall that an $R$-module $M$ is said to be $S$-flat if for any finitely generated ideal $I$ of $R$, the natural homomorphism $(I \otimes_{R} M)_S \to (R \otimes_{R} M)_S$ is a monomorphism \cite[Definition 2.5.]{Qi1}; equivalently, $M_{S}$ is a flat $R_{S}$-module \cite[Proposition 2.6]{Qi1}. Notice that any flat $R$-module is $S$-flat. A general framework for $S$-flat modules was developed in the paper \cite{Bou1}.

Motivated by the work \cite{Qi1}, we aim to define an \(S\)-version of injective modules, thereby extending the following well-known characterizations of Noetherian rings to the broader context of \(S\)-Noetherian rings:

\noindent

\begin{itemize}
    \item \textbf{Cartan–Eilenberg–Bass theorem} \cite[Theorem 4.3.4]{Wan1}: A ring \(R\) is Noetherian if and only if any direct sum of injective \(R\)-modules is injective, or equivalently, if every direct limit of injective \(R\)-modules over a directed set remains injective.
    
    \item \textbf{Cheatham and Stone's theorem} \cite[Theorem 2]{Che1}: A ring \(R\) is Noetherian if and only if, for any \(R\)-module \(M\), \(M\) is injective if and only if \(M^{++}\) is injective, or equivalently, \(M\) is injective if and only if \(M^+\) is flat.
\end{itemize}

 \medskip
To achieve our aim, we introduce an $S$-version of injectivity, which we call $S$-injective modules. Notice that it is different from the notion of $S$-injective modules in the sense of \cite{Bae1}. Indeed, our $S$-injectivity allows us to establish the $S$-version of Baer's Criterion, which plays a crucial role in our study (see Proposition \ref{2-pro-S-baer}). This will be done in Section $2$, which is devoted to investigating the basic properties of the introduced $S$-injectivity. Namely, we provide homological characterizations of $S$-injective modules similar to those of injective modules (see Propositions \ref{2-prop-s-inj-ext}  and \ref{2-prop-s-ijnj-with-exac-s-exact-seq}). We also demonstrate that the class of $S$-injective modules is closed under direct summands and direct products. Section $3$ is devoted to the $S$-version of the Cartan-Eilenberg-Bass theorem (Corollaries \ref{2-cor-s-noe-vs-s-inj} and \ref{2-cor-rs-is-neo-vis-direct-lim-s-inj}) as well as Cheatham and Stone's theorem (Theorem \ref{2-thm-s-noe-vs-chara}).

\section{ Definition and basic properties of $S$-injective modules}
  Let us begin with:
 
\begin{defn}\label{2-def-S-inj}

 An $R$-module $E$ is said to be $S$-injective if, whenever $i: A\to B$ is a monomorphism and
$h: A_S\to E$ is any morphism of $R$-modules, there exists a morphism of $R$-modules  $g$ making the following diagram commute:
$$\xymatrix{  & E & \\
0\ar[r]& A_{S} \ar[r]^{i_{S}}\ar[u]^{h} & B_{S} \ar@{-->}[lu]_{g}}$$

\end{defn}

Obviously, every injective $R$-module $M$ is $S$-injective. Next, in Example \ref{2-exp-S-inj-non-inj}, we provide an example of an $S$-injective module that is not injective. However, these two concepts coincide for $R_S$-modules,  as we will show in Proposition \ref{2-prop-rs-inj=s-inj}. The canonical ring homomorphism $\theta: R \to R_{S}$ makes every $R_{S}$-module an $R$-module via the module action $r.m = \frac{r}{1}.m$, where $r\in R $ and $m\in M$. Recall from \cite[page 417 (2)]{Dad1} that an $R_S$-module is injective as $R_S$-module if and only if it is injective as $R$-module.    

\begin{prop}\label{2-prop-rs-inj=s-inj} An $R_{S}$-module $E$ is injective as an $R$-module if and only if it is $S$-injective. 
\end{prop}

\proof The "only if" part always holds.

Regarding the "if" part, as discussed above, it suffices to show that $E$ is an injective $R_S$-module. But, this is an immediate consequence of \cite[Corollary 4.79]{Rot1}, which states that every $R_S$-module $M$ is naturally isomorphic to its localization $M_S$ as $R_S$-modules. Additionally, we have the fact that:

\begin{center}
$\Hom_{R_{S}}(M, N) = \Hom_{R}(M, N)$ 
\end{center} for all $R_S$-modules $M$ and $N$.\cqfd


It is worth noting that an $R$-module $E$ is injective if and only if every $R$-morphism $f : I \to E$, where $I$ is an ideal of $R$, can be extended to $R$ (Baer's Criterion), \cite[Theorem 1.1.6]{Gla1}. Replacing "injective" with "$S$-injective",  we can prove the following result:

\begin{prop}\label{2-pro-S-baer}
An $R$-module $E$ is $S$-injective if and only if every $R$-morphism $f : I_{S}\to E$, where $I$ is an ideal of $R$, can be extended to $R_{S}$.
\end{prop}
 
 \proof We imitate the proof given by Baer with some adaptations.
The "only if" part is straightforward. \\
 
 For the "if" part, consider the following diagram:
$$\xymatrix{  & E & \\
0\ar[r]& A_{S} \ar[r]^{i_{S}}\ar[u]^{f} & B_{S}}$$ where $A$ is a submodule of an $R$-module $B$ and  $i$ is the inclusion.  Let $\X$ be the set of all ordered pairs $(A', g')$, where $A \subseteq A' \subseteq B$ and $g' : A'_{S} \to E$ extends $f$; i.e.,  $g^{'}|A_{S} = f $. Note that $\X\neq \varnothing$ because $(A,f)\in \X$. Partially order $\X$ by defining \begin{center}
$(A', g') \leq (A'', g'')$
\end{center} to mean $A'\subseteq  A''$ and $g''$ extends $g'$. We may prove easily that chains in $\X$ have upper bounds in $\X$; hence, Zorn's lemma applies, and there exists a maximal element $(M, m)$ in $\X$. If $M_{S} = B_{S} $, we are done, and so we may assume that there is $ b\in B $ with $ \frac{b}{1}\notin M_S $. \\ Define \begin{center}
$I = \{r \in R / rb\in M\}$.
\end{center}
It is easy to see that $I$ is an ideal of $R$. Define $h : I_{S} \to E$ by \begin{center}
$h(\frac{a}{s}) = m(\frac{ab}{s})$.
\end{center}
By hypothesis, there is a map $h^* : R_{S} \to E$ extending $h$. Finally, define $M' = M + <b> $ and $m' : M'_S \to E$ by
\begin{center}
$m'(\frac{a+\alpha b}{s}) = m(\frac{a}{s}) + h^*(\frac{\alpha}{s})$,
\end{center} where $\alpha \in R$, $s\in S$, and $a\in M$.
Let us show that $m'$ is well-defined. If $\frac{a+\alpha b}{s} = \frac{a'+\alpha' b}{s'}$, then \begin{center}
$\frac{a}{s} - \frac{a'}{s'} = \frac{\alpha' b}{s'} - \frac{\alpha b}{s} = \frac{(\alpha' s-\alpha s') b}{ss'} \in M_{S} $,
\end{center} so there exists $n\in M$ and $r\in S$ such that
 $\frac{n}{r} = \frac{(\alpha' s-\alpha s') b}{ss'} \in M_{S} $ and then $lss'n= lr(\alpha' s-\alpha s') b$ for some $l\in S$; it follows that $lr(\alpha' s-\alpha s')  \in I$. Therefore, $h(\frac{lr(\alpha' s-\alpha s') }{lrss'})$ is defined, and we have
 $m(\frac{a}{s}) - m(\frac{a'}{s'}) = m(\frac{lr (\alpha' s-\alpha s') b}{lrss'})=h(\frac{lr (\alpha' s-\alpha s') }{lrss'})= h^*( (\frac{lr (\alpha' s-\alpha s') }{lrss'}))=h^*(\frac{\alpha' }{s'})-h^*(\frac{\alpha }{s})$. Thus, $m(\frac{a}{s})+ h^*(\frac{\alpha }{s}) = m(\frac{a'}{s'})+ h^*(\frac{\alpha' }{s'})$ 
 as desired. Clearly, $(M',m')\in \X$ and $m'(\frac{a}{s}) = m(\frac{a}{s})$ for all $a \in M$ and $s \in S$, so that the map $m'$ extends $m$. We conclude that $(M, m) < (M', m')$, contradicting the maximality of $(M, m)$. Therefore, $M_{S} = B_{S}$, the map $m$ is a lifting of $f$, and $E$ is $S$-injective. 
 \cqfd

\begin{prop}\label{2-prop-s-inj-ext} Let $M$ be an $R$-module.  Consider the following assertions:
\begin{enumerate}
\item $\Ext_R^1(N_S,M)=0$ for any $R$-module $N$. 
\item $ \Ext_R^1(R_S/I_S,M)=0$ for any ideal $I$ of $R$. 
\item $M$ is $S$-injective.
\end{enumerate}

The implications $1.\Rightarrow 2.\Rightarrow 3.$ hold true. Assuming that $R_S$ is projective as an $R$-module, then all the three  assertions are equivalent.

\end{prop}

\proof $1\Rightarrow 2.$ is trivial. 

$2\Rightarrow 3.$ Follows by Proposition \ref{2-pro-S-baer}.

$ 3\Rightarrow 1.$ Assume that $R_S$ is projective. Let $N$ be an $R$-module. There exists an exact sequence of $R_S$-modules: $$0\to K \to P\to N_S\to 0,$$
where $P$ is a projective $R_S$-module. This gives rise to the exact sequence
$$\Hom_R(P,M)\to\Hom_S(K,M)\to \Ext_R^1(N_S,M)\to\Ext_R^1(P,M).$$
Since $R_S$ is projective, $P$ is also projective, and hence $\Ext^1_R(P, M) = 0$. Therefore, $\Ext^1_R(N_S, M) = 0$, because the homomorphism $\Hom_R(P,M)\to\Hom_R(K,M)$ is surjective. \cqfd

  Recall from \cite[Theorem 3.10.22]{Wan1} that a commutative ring $R$ is perfect if and only if every flat $R$-module is projective. Since $R_S$ is a flat $R$-module \cite[Theorem 4.80]{Rot1}, the following result is an immediate consequence of Proposition \ref{2-prop-s-inj-ext}.
  
 \begin{cor} Assume that $R$ is perfect. Then,  an $R$-module $M$ is $S$-injective if and only if $\Ext_R^1(R_S/I_S,M)=0$ for any ideals $I$ of $R$.
 \end{cor}
  
As in the classical case, $S$-injective modules can be characterized through short exact sequences.
 Recall from \cite[Definition 2.1]{Qi1} that a sequence $ 0 \to A \to B \to C \to 0$ of $R$-modules is said to be $S$-exact if the induced sequence $0 \to A_S \to B_S \to  C_S \to 0$ is exact.  
\begin{prop}\label{2-prop-s-ijnj-with-exac-s-exact-seq} The following statements are equivalent for an $R$-module $M.$
\begin{enumerate}
\item $M$ is $S$-injective.
\item  For every  exact sequence of $R$-modules $ 0\to A\to B\to C \to 0$, the induced sequence \begin{center}
$ 0\to \Hom_{R}(C_{S},M)\to \Hom_{R}(B_{S},M)\to \Hom_{R}(A_{S},M) \to 0 $
\end{center} is exact.
\item For every $S$-exact sequence of $R$-modules $ 0\to A\to B\to C \to 0$, the induced  sequence \begin{center}
$ 0\to \Hom_{R}(C_{S},M)\to \Hom_{R}(B_{S},M)\to \Hom_{R}(A_{S},M) \to 0 $
\end{center} is exact.
\end{enumerate}
\end{prop}

\proof  
$1. \Rightarrow 2.$ Assume that $M$ is $S$-injective. Let 
\begin{center}
$\xymatrix@!0 @R=7mm  @C=1.3cm { 0\ar[r]& A \ar[r]^{i} & B \ar[r]^{f}& C\ar[r]^{}&0}$

\end{center}be a short exact sequence. We prove the exactness of:
\begin{center}
$\xymatrix@!0 @R=5mm @C=2.9cm {\hspace{1.7cm} 0 \ar[r]^{} & \Hom_R(C_S, M) \ar[r]^{f_S^*}& \Hom_R(B_S, M) \ar[r]^{i_S^*}& \Hom_R(A_S, M) \ar[r]^{} &0.\hspace{1.7cm}}$\hspace{1cm}
\end{center} Since $\Hom_R(-, M)$ is a left exact contravariant functor, it suffices to show that $i_S^*$ is surjective. Let $f \in  \Hom_R(A_S, M)$. Since $M$ is $S$-injective, there exists $g \in  \Hom_R(B_S, M)$ with $f =gi_S= i_S^{*}(g) $. Hence, $i_S^*$ is surjective.

$2. \Rightarrow 1.$ Let $i:A\to B$ be a monomorphism, and let $f : A_S \to M$. By hypothesis, the induced homomorphism  
 $i^{*} : \Hom_R(B_S, M)\to  \Hom_R(A_S, M)$ is surjective. Then
 there exists $g: B_S \to M$ such that $gi_S = f$, and as a result, the appropriate diagram commutes. Therefore,  $M$ is $S$-injective.

 $2. \Rightarrow 3.$ Let $\xymatrix@!0 @R=7mm  @C=1cm {  0\ar[r]& A \ar[r]^{i} & B \ar[r]^{f}& C\ar[r]^{}&0}$ be an $S$-exact sequence. We need to show that \begin{center}
$\xymatrix {0 \ar[r] & \Hom_R(C_S, M) \ar[r]^{f_S^*}& \Hom_R(B_S, M) \ar[r]^{i_S^*}& \Hom_R(A_S, M) \ar[r] &0}$
\end{center}
is exact. Since $\xymatrix@!0 @R=7mm  @C=1cm {  0\ar[r]& A_S \ar[r]^{i_S} & B_S \ar[r]^{f_S}& C_S\ar[r]^{}&0}$ is an exact sequence and $\Hom(-,M)$ is a left exact contravariant functor, it suffices to show that $i_S^*$ is surjective. Let $h \in  \Hom_R(A_S, M)$. Consider the following exact sequence: \begin{center}
$0\to \ker(i)\to A\to \Im(i)\to 0.$
\end{center} By (2), the induced sequence 
\begin{center}
$ 0\to \Hom_{R}(\Im(i)_{S},M)\to \Hom_{R}(A_{S},M)\to \Hom_{R}(\Ker(i)_{S},M)=0 $

\end{center} is exact. Then, there is $g\in \Hom_{R}(\Im(i)_{S},M)$ such that $h=gi_S.$                  
 
 Now, the inclusion map $k : \Im(i)\to B$ induces the exact sequence \begin{center}
 $0\to \Im(i)\to B \to B/\Im(i)\to 0$.
 \end{center} Again, by (2), the induced sequence
  
\begin{center}
$ 0\to \Hom_{R}((B/\Im(i))_S,M)\to \Hom_{R}(B_{S},M)\to \Hom_{R}(\Im(i)_{S},M)\to 0 $

\end{center}  is exact. So  there exists $g'\in \Hom_{R}(B_{S},M)$ such that $g=g'k_S$. Finally, $h=(g'k_S)i_S=g'(k_Si_S)= g'i_S=i^*_S(g')$, which means that $i_S^*$ is surjective. 

$3. \Rightarrow 2.$ Since $R_S$ is a flat $R$-module, every exact sequence is $S$-exact.  \cqfd

 
 We have the following interesting consequence:
\begin{prop}\label{2-prop-S-inj-colocalization-inj} 
For any \( R \)-module \( M \), \( M \) is \( S \)-injective if and only if \( \Hom_R(R_S, M) \) is injective.
\end{prop}
\proof This follows from Proposition \ref{2-prop-s-ijnj-with-exac-s-exact-seq} and the natural isomorphism \begin{center}
$\Hom_R(A, \Hom_R(B,C)) \cong  \Hom_R(A\otimes_R B,C)$,
\end{center} for any $R$-modules $A$, $B,$ and $C$ \cite[Theorem 2.75]{Rot1}. 
\cqfd

 Notice that, using the natural isomorphism $\Hom_R(A, \Hom_R(B,C)) \cong \Hom_R(A \otimes_R B,C)$, where $A$, $B$, and $C$ are arbitrary $R$-modules (see \cite[Theorem 2.75]{Rot1}), along with Proposition \ref{2-prop-S-inj-colocalization-inj} and Baer's criterion, we obtain another proof of Proposition \ref{2-pro-S-baer}.

\medskip

Proposition \ref{2-prop-S-inj-colocalization-inj} allows us to demonstrate that $S$-injectivity behaves similarly to classical injectivity with respect to direct products.

 \begin{prop}\label{2-prop-prod-of-s-inj}
Let $(M_i)_{i\in I}$ be a family of $R$-modules. Then $\prod\limits_{i\in I} M_i$ is $S$-injective if and only if each $M_i$ is $S$-injective. In particular, every direct summand of an $S$-injective $R$-module is $S$-injective.

\end{prop}
\proof By Proposition \ref{2-prop-S-inj-colocalization-inj}, $\prod\limits_{i\in I} M_i$ is $S$-injective if and only if $\Hom_R(R_S,\prod\limits_{i\in I} M_i)$ is injective. However, since $\Hom_R(R_S,\prod\limits_{i\in I} M_i)\cong \prod\limits_{i\in I} \Hom_R(R_S,M_i)$ by \cite[Theorem 2.30]{Rot1}, it follows from \cite[Proposition 3.28]{Rot1} that $\Hom_R(R_S,\prod\limits_{i\in I} M_i)$ is injective if and only if $\Hom_R(R_S,M_i)$ is injective for each $i\in I$.  Again by Proposition \ref{2-prop-S-inj-colocalization-inj}, this hold if and only if $M_i$ is $S$-injective for any $i\in I$. \cqfd

Recall that an $R$-module $M$ is said to be $S$-flat if for any finitely generated ideal $I$ of $R$, the natural homomorphism $I \otimes_{R} M \to R \otimes_{R} M $ is an $S$-monomorphism; equivalently, $M_{S}$ is a flat $R_{S}$-module \cite[Proposition 2.6]{Qi1}. 

 Recall the Lambek's characterization of flat modules: An $R$-module $M$ is flat if and only if its character $\Hom_{\Z}(M,\Q/\Z)$ is an injective $R$-module \cite[Theorem 1.2.1]{Gla1}. Here, we prove its $S$-version.

\begin{prop}\label{2-prop-lambek}
The following assertions are equivalent for an $R$-module $M$:
\begin{enumerate}
\item $M$ is $S$-flat.
\item $\Hom_{\Z}(M,\Q/\Z)$ is $S$-injective. 
\end{enumerate}
\end{prop}

\proof  This follows from the following natural isomorphisms: 
$$\Hom_{R}(\E_{S}, \Hom_{\Z}(M, \Q/\Z))\cong \Hom_{R}(\E_S\otimes_{R} M, \Q/\Z) \cong \Hom_{R}((\E\otimes_{R} M)_{S}, \Q/\Z),$$ where $\E$
is a short exact sequence of $R$-modules, and the fact that $\E$ is exact if and only if $\Hom_{\Z}(\E,\Q/\Z)$ is exact \cite[Lemma 3.53]{Rot1}. \cqfd

We use Proposition \ref{2-prop-lambek} to give an example of an $S$-injective $R$-module which is not injective:

\begin{exmp}\label{2-exp-S-inj-non-inj} Let $M$ be an $S$-flat module which is not flat \cite{Qi1}. Then , by Proposition \ref{2-prop-lambek}, $\Hom_{\Z}(M, \Q /\Z) $ is an $S$-injective $R$-module, but it is not injective by \cite[Theorem 1.2.1]{Gla1}. 
\end{exmp}


\section{Applications}

In this section, we are interested in the \(S\)-version of the two classical results: the Cartan–Eilenberg–Bass theorem and Cheatham and Stone’s theorem.

It is clear from  Proposition \ref{2-prop-prod-of-s-inj} that  a finite direct sum of $S$-injective $R$-modules is also $S$-injective. We start this section by extending this fact under some conditions.


\begin{prop}\label{2-prop-s-noe-implq-dirct(sum-inj-is-inj}
Assume that $R_S$ is a Noetherian ring and that $R_S$ is finitely generated as an $R$-module. If $(M_i)_{i\in J}$ is a family of $S$-injective $R$-modules, then $\bigoplus\limits_{j \in  J}  M_j $ is an $S$-injective $R$-module.

\end{prop}

\proof We imitate the proof given by \cite[Proposition 3.31]{Rot1}. By Proposition \ref{2-pro-S-baer}, it suffices to complete the diagram
$$\xymatrix{  & \bigoplus\limits_{j\in J} M_j & \\
0\ar[r]& I_{S} \ar[r]^{i_{S}}\ar[u]^{f} & R_{S} \ar@{-->}[lu]_{g}}$$ where $I$ is an ideal of $R$. If $x \in \bigoplus\limits_{j\in J} M_j$, then $x = (e_j)_{j\in J}$, where, for each $j\in  j$,  $e_j \in M_j$. Let $\text{Supp}(x) = \{j\in J : e_j\neq 0\}$.  Since $R_S$ is Noetherian, the ideal $I_S$ is finitely generated as an $R_S$-module. Moreover, since $R_S$ is finitely generated as an $R$-module, $I_S$ is finitely generated as an $R$-module as well. Thus, we can write $I_S = Rx_1 + \cdots + Rx_n$.  Since, for each $ k\in \{1, \ldots, n\}$, \(f(x_k)\) has finite support $\text{Supp}(f(x_k))\subset J$, the set $J' =\bigcup\limits_{k=1}^{k=n} \text{Supp}( f(x_k)) $ is a finite set, and $\Im(f) \subseteq \bigoplus\limits_{j\in J'} M_j$. By Proposition \ref{2-prop-prod-of-s-inj},
this finite direct sum is $S$-injective. Hence, there is an $R$-morphism $g' : R_S\to \bigoplus\limits_{j\in J'} M_j$ extending $f$. Composing $g'$ with the inclusion of $\bigoplus\limits_{j'\in J'} M_j'$  into $\bigoplus\limits_{j\in J} M_j$ completes the given diagram. \cqfd

\begin{cor} Let \( R \) be a ring and \( S = \{s_1, s_2, \ldots, s_n\} \subseteq R \) be a finite multiplicative subset of \( R \). If $R$ is $S$-Noetherian, then every direct sum of $S$-injective $R$-modules is $S$-injective.
\end{cor}

Recall that an injective $R$-module $M$ is said to be $\Sigma$-injective if every direct sum
of copies of $M$ is injective \cite[Definition 4.3.1]{Wan1}. We say that an $S$-injective $R$-module $M$ is $\Sigma$-$S$-injective if every direct sum of copies of $M$ is $S$-injective.  

\begin{prop}\label{2-prop-rs-is-neo-vis-sum-s-inj} Let $S$ be a multiplicative subset of $R$ such that $R_S$ is finitely generated as an $R$-module. The following assertions are equivalent:
\begin{enumerate}
\item $R_S$ is Noetherian.
\item Every direct sum of $S$-injective $R$-modules is $S$-injective.
\item Every direct sum of a countably infinite family of $S$-injective $R$-modules is $S$-injective.
\item Every direct sum of a countably infinite family of injective $R$-modules is $S$-injective.
\item Every direct sum of a countably infinite family of injective $R_S$-modules is $S$-injective.

\item Every $S$-injective $R$-module is $\Sigma$-$S$-injective.
\item Every injective $R$-module is $\Sigma$-$S$-injective.
\item Every injective $R_S$-module is $\Sigma$-$S$-injective.
\end{enumerate}
\end{prop}
\proof $1.\Rightarrow 2.$  This follows from Proposition \ref{2-prop-s-noe-implq-dirct(sum-inj-is-inj}.

$2.\Rightarrow 3.\Rightarrow 4.$ and $2.\Rightarrow 6.\Rightarrow 7.$ are trivial. 

$4. \Rightarrow 5.$ and $7.\Rightarrow 8.$ Follow from the fact that every injective $R_S$-module is injective as an $R$-module \cite[page 417 (2)]{Dad1}.

$5.\Rightarrow 1.$ and $8.\Rightarrow 1.$ Follow from Proposition \ref{2-prop-rs-inj=s-inj}, \cite[Theorem 4.3.4]{Wan1}, and

the fact that an $R_S$-module is injective if and only if it is injective as an $R$-module \cite[page 417 (2)]{Dad1}.\cqfd

Given a commutative ring $R$ and a multiplicative subset $S \subseteq R$, we say that the $S$-torsion in $R$ is bounded by $s_0 \in S$, if for every $s \in S$ and $r \in  R$,  whenever $sr = 0$, it follows that  $s_0r = 0$. This definition can be found in \cite{Pos1}. If $S$ is finite, then the $S$-torsion is bounded by the product of all elements of $S$.

\begin{lem}\label{2-lem-s-noeth-rs-noe}
Let $R$ be a commutative ring and $S \subseteq R$ be a multiplicative subset such that the $S$-torsion in $R$ is bounded by $s_0$. Assume that $R_S$ is finitely generated as $R$-module. Then $R$ is $S$-Noetherian if and only if $R_S$ is Noetherian.
\end{lem}
\proof The "only if" part always holds. To prove the "if" part, suppose that
we are given an ideal $I$ of $R$.  Since $R_S$ is Noetherian and finitely generated as an $R$-module, $I_S=R\frac{a_1}{s_1}+\cdots+ R\frac{a_m}{s_m}$, for some $a_1,..., a_m\in I$, and $s_1,..., s_m\in S$. Let $a\in I$ and $t_0=s_1s_2\cdots s_m$. There  exist $\alpha_1,...,\alpha_m \in R$  such that
 \begin{center}
  $\frac{a}{t_0}=\alpha_1\frac{a_1}{s_1}+\cdots+ \alpha_m\frac{a_m}{s_m}=\frac{\beta_1 a_1+\cdots+\beta_m a_m}{t_0}$
  \end{center} for some $\beta_i\in R$. Then, there exists $t'\in S$ such that $t't_0(a-(\beta_1 a_1+\cdots\beta_m a_m))=0$. Hence, $s_0(a-(\beta_1 a_1+\cdots\beta_m a_m))=0$. Then $s_0a\in I'=Ra_1+\cdots +Ra_m\subseteq I$.  Hence, $I$ is $S$-finite. Therefore, $R$ is $S$-Noetherian. \cqfd

\begin{cor}
Let \( R \) be a ring. For every finite multiplicative subset \( S = \{s_1, \ldots, s_n\} \) of \( R \), \( R \) is \( S \)-Noetherian if and only if \( R_S \) is Noetherian.
\end{cor}

\proof The $S$-torsion in $R$ is bounded by $s_0=s_1s_2\cdots s_n$. Moreover, $R_S=R\frac{1}{s_1}+\cdots +R\frac{1}{s_n}$. Thus, the result follow immediately  from Lemma \ref{2-lem-s-noeth-rs-noe}. \cqfd

We deduce the following result, which may be viewed as an extension of \cite[Theorem 4.3.4]{Wan1}, when the $S$-torsion in $R$ is bounded and $R_S$ is finitely generated. 

\begin{cor}\label{2-cor-s-noe-vs-s-inj} Let $R$ be a commutative ring and $S \subseteq R$ be a multiplicative subset such that the $S$-torsion in $R$ is bounded. Assume that $R_S$ is finitely generated as an $R$-module. Then, the following statements are equivalent:
\begin{enumerate}
\item $R$ is $S$-Noetherian.
\item Every direct sum of $S$-injective $R$-modules is $S$-injective.
\item Every direct sum of countably infinite $S$-injective $R$-modules is $S$-injective.
\item Every $S$-injective $R$-module is $\Sigma$-$S$-injective.
\end{enumerate}
\end{cor}
\proof  This follows by Proposition \ref{2-prop-rs-is-neo-vis-sum-s-inj} and Lemma \ref{2-lem-s-noeth-rs-noe}. \cqfd

\begin{cor}\label{2-cor-rs-is-neo-vis-direct-lim-s-inj} Let $R$ be a commutative  ring and $S \subseteq R$ be a multiplicative subset such that the $S$-torsion in $R$ is bounded by $s_0$. Assume that $R_S$ is finitely presented as $R$-module. Then, the following statements are equivalent:
\begin{enumerate}
\item $R$ is $S$-Noetherian.
\item Every direct limit of $S$-injective $R$-modules over a directed set is $S$-injective.
\end{enumerate}

\end{cor}
\proof $1.\Rightarrow 2.$  Let $(M_i)_{i\in J}$  be a direct system of $S$-injective modules over a directed set $J$. Let $I$ be an ideal of $R$. Since $R$ is $S$-Noetherian, $I$ is $S$-finite. Then, $I_S$ is  finitely generated as $R_S$-module. Since $R_S$ is finitely generated, $I_S$ is finitely generated as an $R$-module.  By \cite[Theorem 2.1.2]{Gla1}, $R_S/I_S$ is a finitely presented $R$-module.   By \cite[ Theorem 3.9.4]{Wan1},

\[\Ext^1_R(R_S/I_S, \varinjlim M_i)\cong  \varinjlim\Ext^1_R(R_S/I_S,  M_i)=0.\] By \cite[Theorem 3.56]{Rot1}, $R_S$ is a projective $R$-module. Therefore, it follows from Proposition \ref{2-prop-s-inj-ext} that $\varinjlim M_i$ is $S$-injective.

$2.\Rightarrow 1.$ By \cite[Example 2.5.30]{Wan1}, every direct sum of $S$-injective modules is a direct limit of $S$-injective modules over a directed set. Hence, $R$ is $S$-Noetherian by Corollary \ref{2-cor-s-noe-vs-s-inj}. \cqfd



We now present an $S$-counterpart of the classical result by Cheatham and Stone \cite[Theorem 2]{Che1}. To establish this, we first prove the following lemma.

\begin{lem}\label{lem-4} Let $R$ be a ring and $S$ a multiplicative subset of $R$ such that $R_S$ is a finitely presented $R$-module. Assume that  $R_S$ is a coherent ring. Then, for any $R$-module $M$, any $S$-finitely presented $R$-module $N$, and any $n \geq 0$:   $$ \Tor_R^n(M^+, N_S) \cong \Ext_R^n(N_S,M)^+.$$
\end{lem} 
 
\proof 
Let $N$ be an $S$-finitely presented $R$-module. Then, $N_S$ is a finitely presented $R_S$-module by \cite[Remark 3.4]{Ben1}. Since $R_S$ is coherent, $N_S$ has a projective resolution composed of finitely generated $R_S$-modules \cite[Corollary 2.5.2]{Gla1}. On the other hand, as $R_S$ is a finitely generated projective $R$-module \cite[Theorem 3.56]{Rot1}, every finitely generated projective $R_S$-module is also a finitely generated projective $R$-module. Therefore, $N_S$ has a projective resolution composed of finitely generated $R$-modules. Consequently, the result follows from \cite[Theorem 1.1.8]{Gla1}. \cqfd

 \begin{thm}\label{2-thm-s-noe-vs-chara} Let $R$ be a commutative ring and $S \subseteq R$ be a multiplicative subset such that the $S$-torsion in $R$ is bounded. Assume that $R_S$ is finitely presented as an $R$-module. Then, the following statements are equivalent:
 
\begin{enumerate}
\item $R$ is $S$-Noetherian.
\item $M$ is $S$-injective if and only if $M^{++}$ is $S$-injective.
 \item $M$ is $S$-injective if and only if $M^{+}$ is $S$-flat.
\end{enumerate}
\end{thm}
\proof  $1.\Rightarrow 3.$ For any ideal $I$ of $R$, there exists a finitely generated subideal $I'$ of $I$ such that $sI\subseteq I'$ for some $s\in S$. Then, $(R/I)_S\cong (R/I')_S$. Consequently, according to Lemma \ref{lem-4},
$$ \Tor_R^1(M^+, (R/I)_S) \cong \Ext_R^1((R/I)_S,M)^+.$$
This holds true for any ideal $I$ of $R$.  Therefore, (3) follows from  Proposition \ref{2-prop-s-inj-ext} and \cite[Proposition 2.5]{Bou2}.

$2.\Leftrightarrow 3.$ Follows from Proposition \ref{2-prop-lambek}.

$3.\Rightarrow 1.$ Using Proposition \ref{2-prop-lambek} and \cite[Theorem 3.6(4)]{Bou4}, one can easily see that if (3) holds, then $R$ is $S$-coherent. Let $(M_i)_{i\in I}$ be a family of $S$-injective $R$-modules. By $(3)$, $M_i^+$ is $S$-flat for any $i\in I$. Since $R$ is $S$-coherent, $\prod\limits_{i\in I} M_i^+$ is $S$-flat by \cite[Theorem 4.2 and Proposition 4.4]{Bou2}. By Proposition \ref{2-prop-lambek}, $(\prod\limits_{i\in I} M_i^+)^+$ is $S$-injective. Since, \begin{center}
$(\prod\limits_{i\in I} M_i^+)^+\cong (\bigoplus\limits_{i\in I}M_i)^{++}$
\end{center}
$\bigoplus\limits_{i\in I}M_i$ is $S$-injective by $(2)$. Therefore, $R$ is $S$-Noetherian by Corollary \ref{2-cor-s-noe-vs-s-inj}. \cqfd 

We conclude this paper with the following example:
 \begin{exmp}
Let $R_1$ be an $S_1$-perfect Noetherian ring $($semisimple ring as an example$)$, $R_2$ be a  commutative ring which is not  Noetherian. Consider the ring $R=R_1\times R_2$ with the multiplicative subset $S = S_1\times 0$. Then 

\begin{enumerate}
\item $R_S\cong (R_1)_{S_1}\times 0$ is a finitely presented projective $R$-module.
\item The $S$-torsion in $R$ is bounded.
\item $R$ is an $S$-Noetherian ring, but it is not Noetherian.
\end{enumerate}  
\end{exmp}
 \proof $1.$ Since $R_1$ is $S_1$-perfect, $(R_1)_{S_1}$ is a finitely generated projective $R_1$-module by \cite[Theorem 4.9]{Bou1}. Then,  $R_S\cong (R_1)_{S_1}\times 0$ is a finitely generated projective $R$-module, so, it is finitely presented.
 
 $2.$ In a commutative Noetherian ring $R$, for any multiplicative subset $S$ of $R$, the $S$-torsion in $R$ is necessarily bounded (see \cite[Page 38]{Pos1}). Thus, the $S_1$-torsion in $R_1$ is bounded by some $s_1\in S_1$. It follows that the $S$-torsion in $R$ is bounded by $(s_1,0)$.
 
$3.$ Obvious. \cqfd

Driss Bennis:  Faculty of Sciences, Mohammed V University in Rabat, Rabat, Morocco.

\noindent e-mail address: driss.bennis@fsr.um5.ac.ma; driss$\_$bennis@hotmail.com

Ayoub Bouziri: Faculty of Sciences, Mohammed V University in Rabat, Rabat, Morocco.

\noindent e-mail address: ayoub$\_$bouziri@um5.ac.ma

\end{document}